\documentclass[a4paper]{amsart}

\usepackage{amsmath,amsthm,amsfonts,amssymb,amscd,xy}

\input xy

\xyoption{all}

\newtheorem{theorem}{Theorem}[section]

\newtheorem{proposition}[theorem]{Proposition}

\newtheorem{lemma}[theorem]{Lemma}

\newtheorem{corollary}[theorem]{Corollary}

\newtheorem{definition}[theorem]{Definition}

\newtheorem{remark}[theorem]{Remark}

\theoremstyle{definition}

\newtheorem{nota}[theorem]{}

\newenvironment{sis}{\left\{\begin{aligned}}{\end{aligned}\right.}

\numberwithin{equation}{section}

\renewcommand{\P}{\mathbb{P}}

\renewcommand{\O}{\mathcal{O}}

\newcommand{\C}{\mathcal{C}}

\newcommand{\F}{\mathcal{F}}

\renewcommand{\L}{\mathcal{L}}

\newcommand{\E}{\mathcal{E}}

\newcommand{\GG}{\mathcal{G}}

\newcommand{\Z}{\mathbb{Z}}

\newcommand{\Q}{\mathbb{Q}}

\DeclareMathOperator{\id}{{\rm id}}

\renewcommand{\det}{{\rm det}}

\renewcommand{\Im}{{\rm Im}}

\newcommand{\Ker}{{\rm Ker}}

\newcommand{\Sym}{{\rm Sym}}

\newcommand{\ov}{\overline}

\newcommand{\HH}{\mathcal{H}_{\rm sm}(1, r, d)}

\newcommand{\G}{\mathbb{G}_m}

\newcommand{\Pic}{{\rm Pic}}

\newcommand{\B}{\mathbb{B}(N)}

\newcommand{\Bsm}{\mathbb{B}_{sm}(N)}

\newcommand{\A}{\mathbb{A}(N)}

\newcommand{\Asm}{\mathbb{A}_{sm}(N)}

\newcommand{\AO}{\mathbb{A}(N)\setminus \{0\}}

\newcommand{\BB}{\mathbb{B}}

\renewcommand{\AA}{\mathbb{A}}

\newcommand{\HHe}{\mathcal{H'}_{\rm sm}(1, r, d)}

\newcommand{\HHa}{\widetilde{\mathcal{H}}_{\rm sm}(1, r, d)}

\renewcommand{\H}{\mathcal{H}_g}

\newcommand{\Ho}{\mathbb{H}}

\newcommand{\Mg}{\mathcal{M}_g}

\begin{document}

\title{The Chow ring of the stack of cyclic covers of the projective line}

\author{Damiano Fulghesu}
\author{Filippo Viviani}
\address{IRMA - Universite de Strasbourg et CNRS 7, rue R. Descartes
67084 Strasbourg, Cedex; France}
\address{Dipartimento di Matematica - Universit\`a degli Studi Roma Tre Largo San Leonardo Murialdo, 1 - 00146 Roma; Italy}
\email{fulghesu@math.unistra.fr, viviani@mat.uniroma3.it}

\thanks{The second author was partially supported by FCT-Ci\^encia2008.}



\begin{abstract}
In this paper we compute the integral Chow ring of the stack of smooth uniform cyclic covers of the projective line and we give explicit generators.
\end{abstract}

\maketitle

\section{Introduction}

The study of moduli spaces of curves has been greatly enhanced by the introduction of algebraic stacks in the work of Deligne and Mumford \cite{DM}. The moduli stack $\mathcal M_g$ of curves of genus $g$ and its compactification $\ov{\mathcal M}_g$ via stable curves are of Deligne-Mumford type. More precisely, they admits a stratification into locally closed substacks
that are quotients of a scheme by an algebraic group acting with finite stabilizers. By using this local structure, Mumford laid down the basis of enumerative computations in \cite{Mum}. Unfortunately, very little is known in general about the integral Chow ring $A^*(\mathcal M _g)$. The intersection rings of $\mathcal{M}_{1,1}$ and $\mathcal{M}_2$ are computed in \cite{EGa} and \cite{Vis}. Even with rational coefficients the ring $A^*(\mathcal M _g)$ is known only
up to $g=5$ (see \cite{Fab} and \cite{Iza}).

In this paper, we give a complete answer for the stack $\HH$ of smooth uniform $\mu_r$-cyclic covers of the projective line, whose branch divisor has degree $N=rd$. In particular, when $r=2$, $d \geq 3$, and the characteristic of the base field $k$ is different from 2, we get the closed substack $\H$ of $\mathcal{M}_g$, whose geometric points are the hyperelliptic curves. The main result we use in order to compute the Chow ring $A^*(\HH )$, is the explicit structure of global quotient stack given by Arsie and Vistoli in  \cite{AV}
$$
\HH = [(\text{Sym}^N(V^*) \backslash \Delta)/G],
$$
where $V$ is the standard representation of $\text{GL}_2$ (so that $\Sym^N(V^*)$ is the vector space of homogeneous binary forms of degree $N$),
$\Delta$ is the discriminant (i.e. the closed subset of binary forms with at least one multiple root),
$G=GL_2$ when $d$ is odd, and $G=\mathbb{G}_m\times PGL_2$ when $d$ is even. We are therefore able to apply equivariant intersection theory developed by Edidin-Graham \cite{EGa} and Totaro \cite{Tot}, getting that
$$
A^*\left( \HH \right) = A^*_G(\text{Sym}^N(V^*) \backslash \Delta).
$$
The first step of our computation is to pass to the projectivization (see Section \ref{proj}).
We reduce to the computation of $A^*_G \left( \mathbb P \left( \text{Sym}^N(V^*) \backslash \Delta \right) \right)$ after showing that
$$
A^*_G(\text{Sym}^N(V^*)\backslash \Delta) = A^*_G \left( \mathbb P \left( \text{Sym}^N(V^*) \backslash \Delta \right) \right) / (t- c_1(\mathcal{D})),
$$
where $\mathcal D$ is a one dimensional representation of $G$
and $t$ is the hyperplane class. Afterward we consider the exact sequence
$$
\xymatrix{
A^{G}_*(\P(\Delta))\ar[r]^(.3){i_*}& A_{G}^*(\mathbb P (\text{Sym}^N(V^*)) )\cong A^*_G[t]/p_N(t) \ar@{->>}[r] &
A^*_G \left( \mathbb P \left( \text{Sym}^N(V^*) \backslash \Delta \right) \right)
}$$
to get (see \ref{reduction-projective})
\begin{equation}
A^*(\HH)=\frac{A^*_G[t]}{\langle t-c_1(\mathcal{D}), p_N(c_1(\mathcal{D})),
{\rm Im}(i_*)\rangle}.
\end{equation}
where $p_N(t)$ is a monic polynomial in the hyperplane class $t$ of degree $N+1$ with coefficients in $A^*_G$. Then we have two cases.
\begin{enumerate}
\item When $d$ is odd, we follow \cite{EF} to show that $(\Im(i_*))$ is generated by two elements and that the class $p_N(c_1(\mathcal{D}))$ belongs to the ideal $(\Im (i_*), t- c_1(\mathcal{D}))$.
If ${\rm char}(k)=0$ or ${\rm char}(k)>N$, we get (Theorem \ref{main-odd})
$$
A^*(\HH)=\frac{\Z[c_1,c_2]}{\langle r(N-1) c_1, \frac{(N-r)(N-r-2)}{4} c_1^2-N(N-2) c_2\rangle}.
$$
where $c_1$ and $c_2$ are (the pull-back of) the Chern classes of the standard representation $V$ of $GL_2$.
\item When $d$ is even, the group $G$ is $\mathbb G_m \times PGL_2$ and $A^*_G$ is $\Z[c_1, c_2, c_3]/<2c_3>$. Since $\mathbb G_m$ acts trivially on $\P( \text{Sym}^N(V^*))$, we can consider only the action of $PGL_2$. In Proposition \ref{Image-s}, we show that $(\Im (i_*))$ has two generators  (which we compute in Proposition \ref{first-map}). It is enough to prove the statement first with coefficients in $\mathbb Z [1/2]$ (\ref{invert2}), then in the localization $\mathbb Z_{(2)}$ (\ref{loc2}).
The last step (see Section \ref{empty}) is to check if the class $p_N(c_1(\mathcal D))$ belongs to the ideal $I=(\Im (i_*), t- c_1(\mathcal{D}))$. The answer depends on $r$. More precisely $p_N(c_1(\mathcal D)) \in I$ if and only if $r$ is even. We compute the ring $A^*(\HH)$ in Theorem \ref{main-even}. If ${\rm char}(k)=0$ or ${\rm char}(k)>N$, we have
$$
A^*(\HH)=\frac{\Z[c_1, c_2, c_3]}{\left\langle p_N(-rc_1), 2c_3, 2r(N-1)c_1, 2r^2 c_1^2-\frac{N(N-2)}{2}c_2 \right\rangle}
$$
from where we can remove $p_N(-rc_1)$ if $r$ is even.
\end{enumerate}

\begin{remark}
In particular, we get that the Picard group of $\HH$ is cyclic of order $r(N-1)$ if $d$ is odd, and of order $2r(N-1)$ if $d$ is even. Therefore, we recover the result of \cite[Thm. 5.1]{AV}.
\end{remark}

In Section \ref{explicit} we give an explicit description for the generators of the Chow ring $A^*(\HH)$ as Chern classes of natural vector bundles (see Theorem \ref{gen-odd} and \ref{gen-even}). Moreover (Proposition \ref{lambda}), we express the tautological $\lambda$  classes of Mumford (\cite{Mum}) in terms of the explicit generators. In the hyperelliptic case we recover the result of \cite{GV}.

Recently, Bolognesi and Vistoli (\cite{BV}) have described the stack $\mathcal{T}_g$
of trigonal curves of genus $g$ as a quotient stack
and have used this explicit presentation to compute the integral Picard group of $\mathcal{T}_g$.
It is likely that some of techniques of this note could be adapted to compute the integral Chow ring
of $\mathcal{T}_g$.

\section{Notations}

Throughout this paper we fix two positive integers $r$ and $d$, and we let $N=rd$. We work over a field $k$ in which $r$ is invertible.

Let us review the description of the stack $\HH$ of smooth uniform $\mu_r$-cyclic covers of the projective line with branch divisor of degree $N=rd$, following Arsie and Vistoli
(see \cite{AV}). We need first to recall the definition of a smooth uniform cyclic cover.

\begin{definition}\cite[Def. 2.1, 2.4]{AV}
\noindent
\begin{enumerate}
\item A uniform $\mu_r$-cyclic cover of a scheme $Y$ consists of a morphism of schemes $f: X \to Y$ together with an action of the group scheme $\mu_r$ on $X$,
such that for each point $q\in Y$, there is an affine neighborhood $V={\rm Spec}R$ of $q\in Y$, together with an element $h\in R$ that is not a zero divisor, and an isomorphism of
$V$-schemes $f^{-1}(V )\cong {\rm Spec} R[x]/(x^r-h)$ which is $\mu_r$-equivariant, when the right hand side is given the obvious actions.  The branch divisor $\Delta_f$ of $f$
is the Cartier divisor on $Y$ whose restriction to ${\rm Spec}(R)$ has local equation $\{h=0\}$.
\item A uniform $\mu_r$-cyclic cover $f:X\to Y$ is said to be smooth if $Y$ and $\Delta_f$ are smooth or, equivalently, if $Y$ and $X$ are smooth.
\end{enumerate}
\end{definition}
The above definition admits a relative version.
\begin{definition}\cite[Def. 2.3, 2.4]{AV}
\noindent
\begin{enumerate}
\item Let $Y\to S$ be a morphism of schemes. A relative uniform $\mu_r$-cyclic cover of $Y\to S$ is a uniform $\mu_r$-cyclic cover $f:X \to Y$ such that the branch divisor $\Delta_f$
is flat over $S$.
\item A relative uniform $\mu_r$-cyclic cover $f:X\to Y$ of $Y\to S$ is said to be smooth if $Y$ and $\Delta_f$ are smooth over $S$ or, equivalently, if $Y$ and $X$ are smooth over $S$.
\end{enumerate}
\end{definition}
Finally, we need to recall the definition of a Brauer-Severi scheme.
\begin{definition}
Let $S$ be a scheme. A Brauer-Severi scheme of relative dimension $n$ over $S$ is a smooth morphism $P\to S$ whose geometric fibers are isomorphic
to the projective space of dimension $n$.
\end{definition}

We can now define $\HH$ as a category fibered in groupoids over the category of $k$-schemes.

\begin{definition}
We denote by $\HH$ the category fibered in groupoids over the category of $k$-schemes, defined as follows.

An object of  $\HH(S)$ over $S$ is a smooth relative uniform $\mu_r$-cyclic cover $X\stackrel{f}{\to} P\to S$ over a Brauer-Severi scheme  $P\to S$ of relative dimension one
such that the branch divisor $\Delta_f$ has relative degree $N=rd$ over $S$.

A morphism from $(X\stackrel{f}{\to}P\to S)$ to $(X'\stackrel{f'}{\to}P'\to S')$ is a commutative diagram
$$\xymatrix{
X\ar[r] \ar[d] & P\ar[r]\ar[d] & S\ar[d] \\
X' \ar[r] & P'\ar[r] & S'
}
$$
where both squares are cartesian and the left hand column is $\mu_r$-equivariant.
\end{definition}

In \cite{AV}, the authors provide an explicit description of $\HH$ as a quotient stack.

\begin{theorem}\cite[Thm. 4.1]{AV}\label{quo-stack}
The category  $\HH$ is isomorphic to the quotient stack
$$
[\Asm/(GL_2/\mu_{d})],
$$
where $\Asm$ is the space of degree $N$ smooth (that is, with distinct roots in the algebraic closure of $k$) binary forms. The group of $d$-th roots of unity $\mu_{d}$ is embedded diagonally into $GL_2$ and the action of $GL_2/\mu_{d}$ on $\Asm$ is given by $[A]\cdot f(x)=f(A^{-1}\cdot x)$.
\end{theorem}
In particular, it follows from the above result that $\HH$ is an irreducible smooth Deligne-Mumford stack of finite type over $k$ of dimension $rd-3$. Analyzing the structure of the algebraic group $GL_2/\mu_{d}$, one can rewrite the above isomorphism of stacks as follows.
\begin{lemma}\label{quotientgroup}\cite[Cor. 4.6]{AV}
\hspace{0cm}
\begin{itemize}
\item[(i)] If $d$ is odd, using the isomorphism $GL_2/\mu_{d} \to GL_2$ given by the map $[A]\mapsto {\rm det}(A)^{\frac{d-1}{2}}A$, the stack $\HH$ can be described as
$$
\HH=[\Asm/GL_2],
$$
with the action given by $A\cdot f(x)={\rm det}(A)^{\frac{r(d-1)}{2}}f(A^{-1}x)$.
\item[(ii)] If $d$ is even, using the isomorphism $GL_2/\mu_{d}\to \mathbb{G}_m\times PGL_2$ given by the map $[A]\mapsto ({\rm det}(A)^{\frac{d}{2}},[A])$, the stack $\HH$ can be described as
$$
\HH=[\Asm/(\mathbb{G}_m\times PGL_2)],
$$
with the action given by $(\alpha,[A])\cdot f(x)=\alpha^{-r}{\rm det}(A)^{\frac{rd}{2}} f(A^{-1}x)$.
\end{itemize}
\end{lemma}
For every smooth scheme $X$ over $k$ endowed with an action of an algebraic group $H$,
we will consider the equivariant
Chow ring $A^*_{H}(X)$ as an algebra over the Chow ring
$A^*_H:=A^*_{H}({\rm Spec}(k))$ of the classifying stack $BH$ of
$H$, via pull-back along the structure map
$X \to {\rm Spec}(k)$. We refer to \cite{EGa} for the definitions and the basic properties of
the equivariant Chow rings.

For the remainder of the paper, we set $G:=GL_2/\mu_d$.
In virtue of Theorem \ref{quo-stack} and
the results of \cite[Sec. 4]{EGa},
the Chow ring $A^*(\HH)$ is isomorphic to the $A^*_G$-algebra $A^*_{G}(\Asm)$.
Using Lemma \ref{quotientgroup}, the Chow ring $A^*_{G}:=A^*_{G}({\rm Spec}(k))$ of the classifying space of
$G=GL_2/\mu_d$ is given as follows (see for example \cite{Pan}).
\begin{lemma}\label{Chow-G}
The  Chow ring $A^*_{G}$ of the classifying stack of $G=GL_2/\mu_d$ is equal to
$$
A^*_G=\begin{cases}
\Z[c_1,c_2] & \text{ if } d \text{ is odd, } \\
 \Z[c_1,c_2,c_3]/(2c_3)  & \text{ if } d \text{ is even. } \\
\end{cases}
$$
\end{lemma}
We can describe the classes appearing in the above Lemma as follows. If $d$ is odd, then $c_1, c_2$ are the Chern classes of the standard representation $V$ of $GL_2$. If $d$ is even, then $c_1$ is the Chern class of the natural representation of $\G$ and $c_2, c_3$ are the Chern classes of the representation ${\rm Sym}^2(V)$ of $PGL_2=SL_2/\mu_2$, where $V$ is the two dimensional standard $k$-representation of $SL_2$ (the first Chern class being trivial).

One last piece of notation: given a set $S$ of elements of $A^*_H(X)$ (for some smooth scheme $X$ acted upon by
an algebraic group $H$), we denote by $\langle S\rangle$ the ideal generated by $S$ inside the ring $A^*_H(X)$.

\section{First reductions}\label{proj}

\subsection{Projectivization}

The first step of our proof consists in passing to the projectivization, as in
 \cite{Vis}. In order to do this, we consider the following $G$-equivariant diagram
\begin{equation}\label{Gm-bundle}
\xymatrix{
\Asm \ar@{^{(}->}[r] \ar[d] & \AO \ar[d] \\
\Bsm \ar@{^{(}->}[r] & \B \\
}
\end{equation}
where $\A$ is the vector space of binary form of degree $N$, $\B:=\P(\AO)$ is its projectivization, and $\Bsm:=\P(\Asm)$ is the open subset of smooth forms.
The vertical arrows of the above diagram are principal $\G$-bundles associated to a $G$-equivariant line bundle $\mathcal{D}\otimes \mathcal{O}(-1)$, where $\mathcal{O}(-1)$ is the tautological bundle on $\B$ and $\mathcal{D}$ is a one dimensional representation of $G$ which, using Lemmas \ref{quotientgroup} and \ref{Chow-G},
has first Chern class in $A^*_G$ given by
\begin{equation}\label{c_1(D)}
 c_1(\mathcal{D})=\begin{cases}
\frac{r(d-1)}{2} c_1=\frac{N-r}{2} c_1& \text{ if } d \text{ is odd, } \\
-r c_1& \text{ if } d \text{ is even. } \\
\end{cases}
\end{equation}
\noindent From the above diagram (\ref{Gm-bundle}), we deduce the following exact diagram of $A^*_G$-algebras
\begin{equation}\label{diag-alg}
\xymatrix{
& \langle t -c_1(\mathcal{D})\rangle \ar@{^{(}->}[d] &
\langle t -c_1(\mathcal{D})\rangle \ar@{^{(}->}[d] \\
A^{G}_*(\P(\Delta))\ar[r]^(.3){i_*}\ar[d]& A_{G}^*(\B)\cong A^*_G[t]/p_N(t) \ar@{->>}[r]
\ar@{->>}[d]& A_{G}^*(\Bsm)\ar@{->>}[d] \\
A^{G}_*(\Delta \backslash \{ 0 \} )\ar[r]& A_{G}^*(\AO) \ar@{->>}[r]
& A_{G}^*(\Asm)
}
\end{equation}
where $\Delta=\A\setminus \Asm$ is the discriminant hypersurface of binary forms having at least two coincident roots over the algebraic closure of $k$,
$t$ is the first $G$-equivariant Chern class of $\mathcal{O}_{\B}(1)$ and $p_N(t)$ is a monic polynomial in $t$ of degree $N+1$ with coefficients in $A^*_G$, whose roots are the opposite of the equivariant Chern roots of the $G$-representation $\A$.
\noindent From the above diagram we deduce that
\begin{equation}\label{reduction-projective}
A^*(\HH)=\frac{A^*_G[t]}{\langle t-c_1(\mathcal{D}), p_N(c_1(\mathcal{D})),
{\rm Im}(i_*)\rangle}.
\end{equation}

\subsection{Stratifying the discriminant} In order to compute the image of the $i_*$, we observe (following \cite{Vis}) that the discriminant hypersurface $\Delta$ has a decreasing filtration into closed subsets
$$
\Delta=\Delta_1\supset \Delta_2 \supset \cdots \supset \Delta_{[N/2]} \supset \Delta_{[N/2]+1}=\emptyset,
$$
where $\Delta_s$ is the closed subset of $\B$ corresponding to forms divisible by the square of a polynomial of degree $s$ over some extension of the base field $k$. There is a natural morphism
$$
\pi_s:\mathbb{B}(s)\times \mathbb{B}(N-2s)\to \B,
$$
which sends $([f],[g])$ into $[f^2g]$. The image of $\pi_s$ is, by definition, the closed subset $\Delta_s$. Arguing as in \cite[Lemma 3.2, 3.3]{Vis} (see also \cite[Prop. 11]{EF}), we conclude that

\begin{lemma}\label{discriminant-strata}
If ${\rm char}(k)=0$ or ${\rm char}(k)>N$, then the ideal $\langle \Im(i_*) \rangle$ is generated by the images of the pushforward maps $\pi_{s*}:A^*_{G}(\mathbb{B}(s)\times \mathbb{B}(N-2s)) \to A^*_{G}(\B)$ where $s=1,\cdots, [N/2]$.
\end{lemma}


\section{$d$ odd}

In the case $d$ odd, we have that $G=GL_2$ and therefore we are in the same situation as in \cite[Section 4]{EF}, where they proved the following
\begin{proposition}\cite[Thm. 19]{EF}\label{Image-GL_2}
If $d$ is odd, then
$$
\langle \Im(\pi_{s*})\rangle_{s\geq 1}=\langle \Im(\pi_{1*})\rangle=\langle 2(N-1)t-N(N-1)c_1, t^2-c_1 t-N(N-2)c_2 \rangle,
$$
where, as usual, $t:=c_1^{GL_2}(\O_{\B}(1))\in A^*_{GL_2}(\B)$.
\end{proposition}
\noindent From this, we can deduce the following
\begin{theorem}\label{main-odd}
If $d$ is odd and ${\rm char}(k)=0$ or ${\rm char}(k)>N$, then
$$
A^*(\HH)=\frac{\Z[c_1,c_2]}{\langle r(N-1) c_1, \frac{(N-r)(N-r-2)}{4} c_1^2-N(N-2) c_2\rangle}.
$$
\end{theorem}
\begin{proof}
The assertion will follow combining (\ref{c_1(D)}), (\ref{reduction-projective}), Lemma \ref{discriminant-strata} and Proposition \ref{Image-GL_2} once we prove that
$$
p_N\left(\frac{N-r}{2}c_1\right)\in \left\langle r(N-1) c_1, \frac{(N-r)(N-r-2)}{4} c_1^2-
N(N-2) c_2\right\rangle.
$$
Indeed, let $V$ be the standard representation of $GL_2$ and let $l_1, l_2$ be its Chern roots. We have $c_1=l_1 + l_2$ and $c_2=l_1 l_2$. We can write more explicitly
$$
\B = \P(\text{Sym}^N (V^*))
$$
Now the Chern roots of  $\Sym^N (V^*)$ are
$$
\{ -(N-i) l_1 - i l_2\}_{i=0, \dots, N}
$$
therefore we have
$$
p_N(t)=\prod_{i=0}^{N} \left( t -(N-i) l_1 - i l_2 \right)
$$
We now multiply the factors of $p_N\left( \frac{N-r}{2} c_1 \right)$ corresponding to $i=0,1,N-1, N$ and we get the polynomial
$$
Q=\left(- \frac{c_1^2}{4} (N - r) (N + r) +
   N^2 c_2 \right) \left(- \frac{c_1^2}{4} (N - r -  2) (N + r - 2) + (N - 2)^2 c_2 \right).
$$
Let us write $f= r(N-1) c_1$ and $g=\frac{(N-r)(N-r-2)}{4} c_1^2-N(N-2) c_2$. We now conclude as in  \cite[Lemma 21]{EF}
$$
Q=f^2c_2 + fgc_1 + g^2.
$$
\end{proof}

\section{$d$ even}

In the case $d$ even, we have that $G=\G\times PGL_2$. Since $\G$ acts trivially on the spaces $\BB(r)$, we can (and we will) work with the equivariant cohomology with respect to $PGL_2$.

Observe that $\BB(s)$ is the projectivization of a $PGL_2$-representation, namely $\AA(s):=\Sym^s(V^*)$, if and only if $s$ is even. In this case, the $PGL_2$-equivariant Chow ring of $\BB(s)$ is equal to
$$
A^*_{PGL_2}(\BB(s))=A^*_{PGL_2}[\xi_s]/(p_s(t)),
$$
where $\xi_s$ is the $PGL_2$-equivariant first Chern class of $\O_{\P^1}(1)$ and $p_s$ is a monic polynomial in $\xi_r$ of degree $s+1$ whose roots are the opposite of the Chern roots of the $PGL_2$-representation $\AA(s)$. Finally observe that $N$ is even since $N=rd$ and we are assuming that $d$ is even. In this case, we set $t:=\xi_N$.

\subsection{Computing ${\rm Im}(\pi_{1*})$}

In this subsection, we compute the image of the map $\pi_{1*}$. We need the following

\begin{lemma}\label{P^1}

The equivariant Chow ring $A^*_{PGL_2}(\P^1)$, considered as an algebra over $A^*_{PGL_2}=\Z[c_2, c_3]/(2c_3)$ is equal to
$$
A_{PGL_2}^*(\P^1)=\Z[c_2, c_3, \tau]/(c_3, \tau^2+c_2)=\Z[\tau],
$$
where $\tau$ is the $PGL_2$-equivariant first Chern class of the $PGL_2$-linearized line bundle $\O_{\P^1}(2)$.
\end{lemma}

\begin{proof}
Since $PGL_2$ acts transitively on $\P^1$ with stabilizer group $E:=\mathbb{G}_a\rtimes \G$, we have that
$$
A^*_{PGL_2}(\P^1)\stackrel{\cong}{\longrightarrow} A^*_E.
$$
According to \cite[Prop. 2.7]{Vez}, we have that
$$
A^*_E \stackrel{\cong}{\longrightarrow} A^*_{\G}=\Z[\tau].
$$
Since $\Pic(B\G)=A^1_{\G}=\Z\cdot \tau$ and $\Pic^{PGL_2}(\P^1)=A^1_{PGL_2}(\P^1)=
\Z\cdot c_1^{PGL_2}(\O_{\P^1}(2))$ (see \cite[Pag. 32-33]{GIT}), we can assume that $\tau=c_1^{PGL_2}(\O_{\P^1}(2))$. Finally, in order to understand the structure of $A^*_{PGL_2}$-algebra of
$A^*_{PGL_2}(\P^1)$, we have to determine the natural pull-back map
$$
A^*_{PGL_2}=\Z[c_2, c_3]/(2c_3)\to A^*_{PGL_2}(\P^1)=\Z[\tau].
$$
Clearly  $c_3$ goes to $0$ since $A^*_{PGL_2}(\P^1)$ is torsion-free, while $c_2$ goes to $\alpha \tau^2$, where $\alpha\in \Z$.

In order to determine $\alpha$, we consider the following commutative diagram
$$
\xymatrix{
 A^*_{PGL_2}=\Z[c_2, c_3]/(2c_3) \ar[r] \ar[d]& A^*_{PGL_2}(\P^1)=\Z[\tau]\ar[d]\\
 A^*_{SL_2}=\Z[c_2] \ar@{^{(}->}[r] & A^*_{SL_2}(\P^1)=\Z[c_2, t]/(t^2+c_2),\\
 }
$$
where $t:=c_1^{SL_2}(\O_{\P^1}(1))$. The right vertical map clearly sends $\tau$ into $2t$, while the left vertical map sends $c_2:=c_2^{PGL_2}(\Sym^2(V))\in A^*_{PGL_2}$ into $c_2^{SL_2}(\Sym^2 V)=4 c_2^{SL_2}(V):=4 c_2\in A^*_{SL_2}$, where $V$ is the standard two dimensional representation of $SL_2$. Therefore, following the image of the element $c_2\in A^*_{PGL_2}$ in the above commutative diagram, we get the equality $4 c_2=4  \alpha t^2=-4\alpha c_2$ in $A^*_{SL_2}(\P^1)$, from which we deduce that $\alpha=-1$.
\end{proof}

\begin{proposition}\label{first-map}
The ideal generated by the image of the push-forward map $\pi_{1*}$ is
equal to
$$
\langle \Im(\pi_{1*}) \rangle=\left\langle 2(N-1)t, 2t^2-\frac{N(N-2)}{2} c_2 \right\rangle.
$$
\end{proposition}

\begin{proof}
By definition of the map $\pi_1$, we have that
$$
\pi_1^*(t)=\tau+\xi_{N-2},
$$
where, with an abuse of notation, we denote with $\tau$ and $\xi_{N-2}$ the pull-back to $A^*_{PGL_2}(\BB(1)\times \BB(N-2))$ of the corresponding classes on the two factors. Therefore, by the projection formula together with the fact that $\tau^2=-c_2$ (see the above Lemma \ref{P^1}), we get that $\Im(\pi_{1*})$ is generated, as an ideal of $A^*_{PGL_2}(\BB(N))$, by $\pi_{1*}(1)$ and $\pi_{1*}(\tau)$.

In order to compute the above two elements $\pi_{1*}(1)$ and $\pi_{1*}(\tau)$,
we will adapt the proof of \cite{Vis} for the case $g=2$. Consider the following $PGL_2$-equivariant commutative diagram
$$
\xymatrix{
\mathbb{B}(1)\times \mathbb{B}(1)^{N-2} \ar^{\delta}[rr] \ar^{\sigma}[d] & &
\mathbb{B}(1)^{N} \ar^{\rho}[d]\\
\mathbb{B}(1)\times \mathbb{B}(N-2) \ar^{\pi_{1}}[rr] & &\B
}
$$
where $\rho$ and $\sigma$ are induced by multiplication of the forms and the map $\delta$ sends $([f_1],[f_2],\cdots,[f_{N-1}])$ into $([f_1],[f_1],[f_2],[f_3],\cdots, [f_{N-1}])$. The maps $\sigma$ and $\rho$ are finite and flat of degree $(N-2)!$ and $N!$ respectively, and therefore we have the formulas $\sigma_*\sigma^*=(N-2)!$ and $\rho_*\rho^*=N!$.

Since the elements $\pi_{1*}(1)$ and $\pi_{1*}(\tau)$ have degree $1$ and $2$ and the Chow ring $A^*_{PGL_2}(\B)$ has no torsion in degree $1$ and $2$, we can work with $\Q$-coefficients. Using the formula $\sigma_*\sigma^*=(N-2)!$, we get (for $i=0, 1$)
\begin{equation}\label{Iformula1}
\pi_{1*}(\tau^i)=\frac{1}{(N-2)!}\rho_*\delta_*\sigma^*(\tau^i).
\end{equation}
Call $\tau_j$ the pullback to $\mathbb{B}(1)^{N}$ of the $PGL_2$-equivariant first Chern class of the line bundle $\mathcal{O}_{\BB(1)}(2)$ on the $j$-th factor (for $j=1,\cdots, N$). By definition of the maps $\sigma$ and $\delta$, we have that $\sigma^*(\tau^i)=\delta^*(\tau_1^i)=\delta^*(\tau_2^i)$ for $i=0,1$. Hence, using the fact that $\delta_*(2)=\tau_1+\tau_2$, we obtain
\begin{equation}\label{Iformula2}
\begin{sis}
& \delta_*\sigma^*(2 )=\delta_*(2)=\tau_1+\tau_2, \\
& \delta_*\sigma^*(2\tau)=\delta_*\delta^*(2\tau_1)=\tau_1(\tau_1+\tau_2).
\end{sis}
\end{equation}
From the definition of the map $\rho$ and the Lemma \ref{P^1}, we get
\begin{equation}\label{Iformula3}
\begin{sis}
& \rho^*(1)=1,\\
& \rho^*(2t)=\sum \tau_j,\\
&\rho^*(4t^2)=\left(\sum \tau_j\right)^2= -N c_2+2 \sum_{j<k}\tau_j\tau_k,
\end{sis}
\end{equation}
where the indices appearing in the above summations range from $1$ to $N$. By taking the pushforward of the above equations (\ref{Iformula3}) and using the formula $\rho_*\rho^*=N!$ and the fact that the pushforwards $\rho_*(\tau_j)$ and $\rho_*(\tau_j\tau_k)$ do not depend on the indices $j$ and $k$ (by symmetry of the map $\rho$), we get
\begin{equation}\label{Iformula4}
\begin{sis}
& \rho_*(1)=N! ,\\
& \rho_*(\tau_j)=2(N-1)! \, t,\\
&\rho_*(\tau_j\tau_k)= (N-2)!\, [4 t^2+N c_2].
\end{sis}
\end{equation}
Putting together the formulas (\ref{Iformula1}), (\ref{Iformula2}), (\ref{Iformula4}),
we get the desired conclusion.
\end{proof}

\subsection{Computing $\langle\pi_{s*}\rangle$} The aim of this subsection is to prove the following
\begin{proposition}\label{Image-s}
If $d$ is even, we have that
$$\langle \Im(\pi_{s*})\rangle_{s\geq 2} \subset \langle \Im(\pi_{1*}) \rangle.
$$
\end{proposition}
Note that it is enough to prove the assertion first with coefficients in the ring $\Z[1/2]$ obtained by inverting $2$ and then in the ring $\Z_{(2)}$ obtained by localizing at $2$. We will treat the two cases separately. We begin with some preliminary Lemmas. The first one generalizes \cite[Prop. 2.3]{Vez}.
\begin{lemma}\label{n-torsion}
Let $X$ be a smooth scheme on which $PGL_n$ acts and consider the induced action of
$SL_n$ via the quotient map $SL_n\twoheadrightarrow PGL_n$. Then the kernel of the natural pull-back map $A^*_{PGL_n}(X) \to A^*_{SL_n}(X)$ is of $n$-torsion.
\end{lemma}
\begin{proof} Recall (see \cite[Prop. 19]{EGa}) that given a group $G$ acting on a smooth scheme $X$, the equivariant Chow ring $A^*_G(X)$ is isomorphic to the operational Chow ring $A^*([X/G])$ of the quotient stack $[X/G]$. More precisely, any element $c\in A^i_G(X)$ defines a compatible system of operations $c(S\to [X/G]): A_*(S)\to A_{*-i}(S)$ for any map $S\to [X/G]$. If $S$ is also smooth then such an operation is induced by the cup product with an element of $A^i(S)$ which, by abuse of notation, we denote also by $c(S\to [X/G])\in A^i(S)$.

Let $p:[X/SL_n]\to [X/PGL_n]$ the natural map of quotient stacks and let $p^*:A^*([X/PGL_n])\to A^*([X/SL_n])$ the pull-back map. Fix an integer $i$ and an element $c \in A^i([X/PGL_n])$ such that $p^*c=0$. We want to prove that $0=n\cdot c\in A^i([X/PGL_n])$. It is enough to prove that for any map $\alpha:S\to [X/PGL_n]$ with $S$ smooth, we have $nc(S\stackrel{\alpha}{\to} [X/PGL_n])=0$.

Indeed, let $V$ be a representation of $PGL_n$ and $U \subseteq V$ an open subset on which $PGL_n$ acts freely and whose complement has codimension higher than $i$. Now if $c$ is an assignment which is 0 on smooth varieties then $c$ is 0 on the torsor $X \times U \to (X \times U)/PGL_n$. But by definition (see \cite[Definition-Proposition 1]{EGa})
$$
A^i((X \times U)/PGL_n) =  A^i([X/PGL_n])
$$
so the class must be 0 in $A^i[X/PGL_n]$.

Now, the map $\alpha:S\to [X/PGL_n]$ corresponds to a $PGL_n$-torsor $P\to S$ over $S$ together with a $PGL_n$-equivariant map $P\to X$. Let $f:\ov{P}\to S$ the Severi-Brauer scheme over $S$ corresponding to the above $PGL_n$-torsor. The map $f:\ov{P}\to S$ is smooth with all the geometric fibers isomorphic to $\P^{n-1}$. Consider the following modified push-forward
$$
\begin{aligned}
 f_{\#}: A^i(\ov{P}) & \longrightarrow A^i(S)\\
\alpha & \mapsto f_*(\alpha\cdot c_{n-1}(T_f)),
\end{aligned}
$$
where $T_f$ is the relative tangent bundle of the map $f: \ov{P}\to S$. By the projection formula and the fact that $\deg c_{n-1}(T_{\P^{n-1}})=n$, we get that
\begin{equation*}
f_{\#}\circ f^*=n.
\end{equation*}
Using this formula, the proof will be complete if we show that
\begin{equation}
f^*c(S\stackrel{\alpha}{\to} [X/PGL_n])=0.
\end{equation}
By functoriality, we have that
$$
f^*c(S\stackrel{\alpha}{\to} [X/PGL_n])= c(\ov{P}\stackrel{\alpha':=f\circ \alpha}{\longrightarrow} [X/PGL_n]).
$$
The map $\alpha':\ov{P}\to [X/PGL_n]$ corresponds to the $PGL_n$-torsor $\ov{P}\times_S P\to \ov{P}$ and the $PGL_n$-equivariant map $\ov{P}\times_S P\to P \to X$. By construction, the $PGL_n$-torsor $\ov{P}\times_S P\to \ov{P}$ is trivial and therefore there exists a $SL_n$-torsor $E \to \ov{P}$ such that $E/\mu_n\cong \ov{P}\times_S P$. Moreover, since the group $\mu_n$ acts trivially on $X$, the $PGL_n$-equivariant map $\ov{P}\times_S P \to X$ extends to a $SL_n$-equivariant map $E\to X$. This data defines a morphism $\beta:\ov{P}\to [X/SL_n]$ such that $\alpha'=p \circ \beta$. Therefore, using the hypothesis that $p^*c=0$, we get that
$$
c(\ov{P}\stackrel{\alpha'}{\longrightarrow} [X/PGL_n])= (p^*c)(\ov{P}\stackrel{\beta}{\longrightarrow} [X/SL_n])=0,
$$
and this concludes the proof.
\end{proof}

\begin{lemma}\label{s-odd}
If $s$ is odd then
$$
\langle \Im(\pi_{s*})\rangle\otimes \Z_{(2)} \subset
\langle \Im(\pi_{1*}) \rangle\otimes \Z_{(2)}.
$$
\end{lemma}

\begin{proof}
Consider the following $PGL_2$-equivariant commutative diagram
$$
\xymatrix{
& \BB(1)\times \BB(s-1)\times \BB(N-2s) \ar^{\rho_s\times \id}[dl]
\ar_{\id\times \pi'_{s-1}}[dr]& \\
\BB(s)\times \BB(N-2s)\ar^{\pi_s}[dr]& &\BB(1)\times \BB(N-2)\ar_{\pi_1}[dl] \\
& \B & \\
}
$$
where $\rho_s$ sends $([f], [g])$ to $[fg]$ and $\pi'_{s-1}$ sends $([f], [g])$ to $[f^2 g]$. The map $\rho_s\times \id$ is finite and flat of degree $s$ and therefore
$$
(\rho_s\times \id)_*\circ (\rho_s\times \id)^*=s \cdot \id.
$$
Using this and the commutativity of the above diagram, we get for any class $\alpha\in A^*_{PGL_2}(\BB(s)\times \BB(N-2s))$
$$
s\cdot \pi_{s*}(\alpha) =\pi_{s*}(\rho_s\times \id)_*(\rho_s\times \id)^*(\alpha)=\pi_{1*}(\id\times \pi'_{s-1})_*(\rho_s\times \id)^*(\alpha).
$$
Since $s$ is odd by hypothesis, and therefore invertible in $\Z_{(2)}$, we deduce that $\pi_{s*}(\alpha)\in \langle \Im(\pi_{1*})\rangle \otimes \Z_{(2)}$.
\end{proof}

\begin{lemma}\label{s-even}
If $s$ is even then $\langle \Im(\pi_{s*})\rangle$ is $2$-divisible in $A^*_{PGL_2}(\B)$.
\end{lemma}

\begin{proof}
We first prove the assertion in the case $s$ is maximal, that is $s=N/2$ and $N$ divisible by $4$.

Consider an element $\alpha\in A^*_{PGL_2}(\BB(s))$. Observe that the element $\pi_{s*}(\alpha)$ is $2$-divisible in the ring $A^*_{PGL_2}(\B)=A^*_{PGL_2}[t]/(p_N(t))$  if and only if the element $t\cdot \pi_{s*}(\alpha)$ is $2$-divisible. The map $\pi_s:\BB(s)\to \B$, sending $[f]$ into $[f^2]$, verifies $\pi_{s}^*(t)=2 \xi_s.$ Therefore, by the projection formula applied to the morphism $\pi_{s}$, we get that
$$
t\pi_{s*}(\alpha)=\pi_{s*}(\pi_s^*(t)\cdot\alpha)=\pi_{s*}(2\xi_{s}\cdot\alpha)=2\pi_{s*}(\xi_{s}\alpha),
$$
from which we conclude in the case $s=N/2$. For the general case, observe that the map $\pi_s$ factors as
$$
\pi_s:\BB(s)\times \BB(N-2s)\stackrel{\pi'_{s}\times \id}{\longrightarrow}
\BB(2s)\times \BB(N-2s)\stackrel{m_s}{\longrightarrow} \B
$$
where $\pi'_{s}$ sends $[f]$ into $[f^2]$ and $m_s$ sends $([f], [g])$ into $[fg]$. It is therefore enough to prove that $\langle \Im (\pi'_s\times \id)_{*}  \rangle $ is $2$-divisible in $A^*_{PGL_2}(\BB(s)\times \BB(N-2s))$.

Observe that since $s$ is even we have an isomorphism
$$
A^*_{PGL_2}(\BB(s)\times \BB(N-2s))=A^*_{PGL_2}(\BB(N-2s))\otimes_{A^*_{PGL_2}} A^*_{PGL_2}(\BB(s)),
$$
and similarly for $A^*_{PGL_2}(\BB(2s)\times \BB(N-2s))$. In particular, the ring $ A^*_{PGL_2}(\BB(s)\times \BB(N-2s))$ (resp. $A^*_{PGL_2}(\BB(2s)\times \BB(N-2s))$) is a $A^*_{PGL_2}(\BB(N-2s))$-module generated by the pull-back along the first projection of the class $\xi_s$ (resp. $\xi_{2s}$). Moreover the push-forward $(\pi'_s\times \id)_*$ is a morphism of $A^*_{PGL_2}(\BB(N-2s))$-modules. Therefore we deduce the $2$-divisibility of the image of the push-forward $(\pi'_s\times \id)_*$ from the previous maximal case.
\end{proof}

We are now ready to prove the Proposition \ref{Image-s}.

\begin{nota}{\emph{\bf Proof of Proposition \ref{Image-s} with coefficients in $\Z[1/2]$}} \label{invert2}

Consider the following commutative diagram
$$
\xymatrix{
A^*_{PGL_2}(\BB(s)\times \BB(N-2s)) \ar[r]^(0,6){\pi_{s*}} \ar[d] &
A^*_{PGL_2}(\BB(N))\ar[d]\\
A^*_{SL_2}(\BB(s)\times \BB(N-2s)) \ar[r]^(0.6){\pi^{SL_2}_{s*}}  &
A^*_{SL_2}(\BB(N)).
}
$$

According to Lemma \ref{n-torsion}, the two vertical arrows become injective after tensoring with $\Z[1/2]$. Therefore in order to prove the inclusion of Proposition \ref{Image-s} with $\Z[1/2]$-coefficients, it is enough to prove the analogous inclusion for the $SL_2$-equivariant Chow rings. But this is proved exactly as in the case $GL_2$ (see \cite[Section 4]{EF}): the same proof works by simply putting $c_1=0$.
\qed
\end{nota}

\begin{nota}{\emph{\bf Proof of Proposition \ref{Image-s} with coefficients in $\Z_{(2)}$}} \label{loc2}

The inclusion $\langle \Im (\pi_{s*})\rangle \subset \langle \Im (\pi_{1*})\rangle$ for $s$ odd follows directly from Lemma \ref{s-odd}. Let us fix an even $s\geq 2$. According to Lemma \ref{s-even}, we have that $\langle \Im(\pi_{s*})\rangle \subset \langle 2\rangle\subset A^*_{PGL_2}(\BB(N))$. A direct check  using Proposition \ref{first-map} and the fact that $N$ is even shows that we have also the inclusion $\langle \Im(\pi_{1*})\rangle \subset \langle 2\rangle$.

Consider the natural pull-back map $p^*:A^*_{PGL_2}(\BB(N))\to A^*_{SL_2}(\BB(N))$. Since $N$ is even, $\BB(N)$ is a projective bundle over the classifying stack $BPGL_2$ and therefore the kernel of $p^*$ is generated by the kernel of the natural pull-back map $A^*_{PGL_2}\to A^*_{SL_2}$, that is $\Ker (p^*)=\langle c_3 \rangle$.

Using the relation $2c_3=0$ it is easy to check that $\langle 2\rangle\cap \langle c_3\rangle=0$. Therefore the pull-back map $p^*$ is injective on the ideal $\langle 2\rangle$. Since $\langle \Im(\pi_{s*})\rangle$ and $\langle \Im(\pi_{1*})\rangle$ are contained in the ideal $\langle 2\rangle$, in order to prove the inclusion $\langle \Im(\pi_{s*})\rangle \subset \langle \Im(\pi_{1*})\rangle$ it is enough to prove the similar inclusion in the ring $A^*_{SL_2}(\BB(N))$. This is done exactly as in the case $GL_2$ (see \cite[Section 4]{EF}): the same proof works by simply putting $c_1=0$.
\qed
\end{nota}

\subsection{Conclusion}
Now we put everything together to prove the following
\begin{theorem}\label{main-even} If $d$ is even and ${\rm char}(k)=0$ or ${\rm char}(k)>N$, we have if $r$ is even

$$
A^*(\HH)=\frac{\Z[c_1, c_2, c_3]}{\left\langle 2c_3, 2r(N-1)c_1, 2r^2 c_1^2-\frac{N(N-2)}{2}c_2 \right\rangle}.
$$
while, if $r$ is odd
$$
A^*(\HH)=\frac{\Z[c_1, c_2, c_3]}{\left\langle p_N(-rc_1), 2c_3, 2r(N-1)c_1, 2r^2 c_1^2-\frac{N(N-2)}{2}c_2 \right\rangle}.
$$
\end{theorem}

\begin{proof}
The assertion will follow combining (\ref{c_1(D)}), (\ref{reduction-projective}), Lemma \ref{discriminant-strata}, Propositions \ref{first-map} and \ref{Image-s} once we prove that
\begin{enumerate}
\item if $r$ is even (Proposition \ref{r-even})
$$
p_N\left(-r c_1 \right) \in \left\langle 2c_3, 2r(N-1)c_1, 2r^2 c_1^2-\frac{N(N-2)}{2}c_2 \right\rangle.
$$
\item if $r$ is odd (Proposition \ref{explicit-r-odd})
$$
p_N\left(-r c_1 \right) \notin \left\langle 2c_3, 2r(N-1)c_1, 2r^2 c_1^2-\frac{N(N-2)}{2}c_2 \right\rangle.
$$
\end{enumerate}
An explicit form of $p_N(-rc_1)$ when $r$ is odd is given in Proposition \ref{explicit-r-odd}.
\end{proof}

\begin{remark}
In the hyperelliptic case $\H=\mathcal{H}_{\rm sm}(1,2,g+1)$, we get the following answer. If $g$ is even (see also \cite{EF})
$$
A^*(\H)=\frac{\Z[c_1,c_2]}{\langle 2(2g+1) c_1, g(g-1) c_1^2-4g(g+1) c_2\rangle}.
$$
If $g$ is odd
$$
A^*(\H)=\frac{\Z[c_1, c_2, c_3]}{\left\langle 2c_3, 4(2g+1)c_1, 8 c_1^2-2g(g+1)c_2 \right\rangle}.
$$
\end{remark}

\section{The polynomial $p_N(t)$}\label{empty}

Throughout this section, we assume that $d$ is even and we set $d=2s$,
so that $N=2rs$.
We want to compute the polynomial $p_N(t)$ in the ring
$$
A^*_{PGL_2}[t]=\Z[c_2,c_3,t]/(2c_3).
$$
Let $V$ be the defining representation of $GL_2$. Let $a,b,c$ be
the Chern roots of $\text{Sym}^2V^*$, seen as a representation of $PGL_2$. We have
$$\begin{sis}
&a + b + c=0,\\
&ab + ac + bc = c_2,\\
&abc=c_3,\\
&2abc=0.
\end{sis}$$
In the next Lemma, we determine the polynomial $p_N(t)$ modulo the ideal $(2)$
of $\Z[c_2,c_3]/(2c_3)$.
\begin{lemma}\label{empty-mod-2}
The polynomial $p_N(t)$ is equivalent modulo $(2)$ to
$$
p_N(t)\equiv
\begin{cases}
 t^{(N+4)/4}(t^3+c_2t+c_3)^{N/4}  \hspace{0,6cm}\mod (2) & \text{ if } \; N \equiv 0 \mod 4 ,\\
 t^{(N-2)/4}(t^3+c_2t+c_3)^{(N+2)/4} \mod (2)& \text{ if } \; N \equiv 2 \mod 4 . \\
 \end{cases}
$$
\end{lemma}
\begin{proof}
By the well-known plethysm  formulas for $\mathrm{PGL}_2$, we have that
$$
\Sym^{2sr}(V^*) \oplus \Sym^{sr - 2}(\Sym^{2}(V^*)) = \Sym^{sr}(\Sym^2(V^*)).
$$
Therefore, from the Whitney's formula and the usual formulas for the Chern roots
of a symmetric product of a representation, we get that
\begin{equation}\label{p-frac}
p_{N}(t)=\frac{\prod_{i,j,k \geq 0}^{i+j+k=rs}(t+ia+jb+kc)}
{\prod_{i,j,k \geq 0}^{i+j+k=rs-2}(t+ia+jb+kc)}.
\end{equation}
Consider the expression of $p_N(t)$ obtained in Lemma \ref{p-frac}.
Since, for all $i,j,k$ such that $i+j+k=rs-2$ we have
$$
(t + i a +j b + k c) \equiv (t + i a + j b + (k+2)c)  \mod (2),
$$
we can simplify the fraction in (\ref{p-frac}) modulo $(2)$ and thus we get
$$
p_N(t) \equiv \prod_{i=0}^{rs}(t + i a + (rs-i)b)
\prod _{i=0}^{rs-1}(t + i a + (rs-1-i)b + c) \mod (2).
$$
We compute separately the two products mod $(2)$. In the first product, the coefficients of $a$ and $b$
have the same parity if $rs$ is even and opposite parity if $rs$ is odd, so that
we get (using that $c\equiv a+b \mod (2)$):
\begin{equation*}
\prod_{i=0}^{rs}(t + i a + (rs-i)b) \equiv
\begin{cases}
t^{\frac{rs}{2}+1}(t + c)^{\frac{rs}{2}} \hspace{1,1cm} \mod (2) &\text{ if } rs \text{ is even, }\\
(t+a)^{\frac{rs+1}{2}}(t+b)^{\frac{rs+1}{2}}  \mod (2) &\text{ if } rs \text{ is odd. }\\
\end{cases}\tag{*}
\end{equation*}
A similar computation for the second product gives
\begin{equation*}
\prod_{i=0}^{rs-1}(t + i a + (rs-1-i)b+ c) \equiv
\begin{cases}
(t+b)^{\frac{rs}{2}}(t + a)^{\frac{rs}{2}} \mod (2) &\text{ if } rs \text{ is even, }\\
t^{\frac{rs-1}{2}}(t+c)^{\frac{rs+1}{2}} \hspace{0,3cm} \mod (2) &\text{ if } rs \text{ is odd. }\\
\end{cases}\tag{**}
\end{equation*}
We conclude by putting together (*) and (**), and using that $N=2rs$ and
$(t+a)(t+b)(t+c)=t^3+c_2t+c_3.$
\end{proof}
We now determine the polynomial $p_N(t)$ modulo the ideal $(c_3)$ of $\Z[c_2,c_3]/(2c_3)$.
\begin{lemma}\label{mult}
The polynomial $p_N(t)$ is equivalent modulo $(c_3)$ to
$$p_N(t)\equiv t \prod_{k=1}^{N/2} (t^2+k^2 c_2) \mod (c_3).
$$
\end{lemma}
\begin{proof}
We compare the $\mathrm{PGL}_2$-equivariant Chow ring $A_{PGL_2}(\B)$ with the
$\mathrm{SL}_2$-equivariant Chow ring $A_{SL_2}(\B)$, in a similar way as
we did in the proof of Lemma \ref{P^1}.
To this aim, let us first compute the Chow ring $A_{SL_2}(\B)$.
Clearly, we have that
$$
A^*_{SL_2}(\B)=A_{SL_2}[\tau]/(q_N(\tau))=\Z[c_2, \tau]/(q_N(\tau)),
$$
where $\tau=c_1^{SL_2}(\O_{\B}(1))$ and $q_N(\tau)$ is a monic polynomial of degree $N+1$ in
$\tau$ with coefficients in $\Z[c_2]$. Let $\alpha$ and $\beta$ be the Chern roots of
$V^*$, seen as a representation of $\mathrm{SL}_2$. We have that
$$\begin{sis}
& \alpha+\beta=0,\\
&\alpha\beta=c_2\in A^*_{SL_2}.
\end{sis}$$
The Chern roots of the $\mathrm{SL}_2$-representation $\Sym^N(V^*)$ are
$\{i\alpha+(N-i)\beta\}_{i=0,\cdots,N}$ and therefore we compute
$$q_N(\tau)=\prod_{i=0}^{N}\left(\tau+i\alpha+(N-i)\beta \right)=$$
$$=\prod_{i=0}^{N/2-1}\left[\left(\tau+i\alpha+(N-i)\beta \right)
\left(\tau+(N-i)\alpha+i\beta \right)\right]\cdot \left(\tau+\frac{N}{2}\alpha+
\frac{N}{2}\beta\right)=$$
$$=\tau \prod_{i=0}^{N/2-1} \left[\tau^2+\left(\frac{N}{2}-i\right)^2 4c_2\right]=
\tau \prod_{k=1}^{N/2} \left[\tau^2+k^24c_2\right].
$$
Now consider the natural commutative diagram of rings (similar to the one
considered in Lemma \ref{P^1}):
$$\xymatrix{
 A^*_{PGL_2}=\Z[c_2, c_3]/(2c_3) \ar[r] \ar[d]& A^*_{PGL_2}(\B)=\Z[c_2,c_3,t]/(2c_3,p_N(t))\ar[d]\\
 A^*_{SL_2}=\Z[c_2] \ar@{^{(}->}[r] & A^*_{SL_2}(\B)=\Z[c_2, \tau]/(q_N(\tau)).\\
}$$
The left vertical maps sends $c_3$ to $0$
and $c_2$ to $4c_2$ (see the proof of Lemma \ref{P^1}), while the right vertical map
obviously sends  $t=c_1^{PGL_2}(\O_{\B}(1))$ into $\tau=c_1^{SL_2}(\O_{\B}(1))$.
This diagram tells us that the polynomial obtained from $p_N(t)$
by substituting $t$ with $\tau$, $c_2$ with $4c_2$ and $c_3$ with $0$
should be equal to $q_N(\tau)$. From the above formula for $q_N(\tau)$,
we get the conclusion.
\end{proof}

We can now put together the previous two Lemmas to get the  following expression
for $p_N(t)\in \Z[c_2,c_3,t]/(2c_3)$.

\begin{corollary}\label{pN-explicit}
If $N\equiv 0 \mod 4$ then we have
$$p_N(t)=t\prod_{k=1}^{N/2} (t^2+k^2 c_2) + t^{\frac{N}{4}+1} \sum_{k=1}^{N/4} \binom{\frac{N}{4}}{k}
(t^3+c_2t)^{\frac{N}{4}-k}c_3^k,$$
while if $N\equiv 2 \mod 4$ then
$$p_N(t)=t\prod_{k=1}^{N/2} (t^2+k^2 c_2) + t^{\frac{N-2}{4}} \sum_{k=1}^{(N+2)/4} \binom{\frac{N+2}{4}}{k}
(t^3+c_2t)^{\frac{N+2}{4}-k}c_3^k.$$
\end{corollary}
\begin{proof}
The polynomial $p_N(t)$ modulo $(c_3)$ is given by Lemma \ref{mult}.
Since $2c_3=0$, the terms of the polynomial that are multiples of $c_3$ are given the
corresponding terms in the expression of $p_N(t)$ modulo $(2)$ (see Lemma \ref{empty-mod-2}).
A straightforward computation allows to conclude.
\end{proof}
Now, we evaluate the class of $p_N(-r c_1)$ modulo the ideal
$$I:=\left\langle 2r(N-1)c_1, 2r^2 c_1^2-\frac{N(N-2)}{2}c_2, 2c_3 \right\rangle
\subset \Z[c_1,c_2,c_3].$$
The answer will depend upon the parity of $r$. Consider first the case where
$r$ is even.
\begin{proposition}\label{r-even}
If $r$ is even then $p_{N}(-r c_1)$ belongs to $I$.
\end{proposition}
\begin{proof}
Substituting $t=-r c_1$ into the first expression of Corollary \ref{pN-explicit}
(note that $N=2dr\equiv 0 \mod 4$) and using the fact that $2c_3=0$ and $r$ is even,
we get that
$$p_N(-rc_1)\equiv -rc_1 \prod_{k=1}^{N/2} (r^2c_1^2+k^2 c_2) \mod I.$$
Consider the element
$$
f:=-r c_1\left(r^2 c^2_1 + \frac{N^2}{4} c_2\right),
$$
which appears as a factor of the above expression for $p_N(-r c_1)$.
We will show that $f\in I$, which will conclude the proof.
Since $r$ is even (and therefore $N=2rs\equiv 0 \mod4$), we can write
$$f=-\frac{r}{2}c_1\left(2 r^2 c^2_1 + \frac{N^2}{2} c_2\right)\equiv
-\frac{r}{2}c_1\left(\frac{N(N-2)}{2}c_2 + \frac{N^2}{2} c_2\right)=
$$
$$=-\frac{N}{4}c_2\cdot 2r(N-1) c_1\equiv 0 \mod I.
$$
\end{proof}

Finally, we consider the case where $r$ is odd.
\begin{proposition}\label{explicit-r-odd}
If $r$ is odd then the expression of $p_N(-rc_1)$ modulo $I$ is equal to
$$p_N(-rc_1)\equiv c_1^{\frac{N}{4}+1} (c_1^3+c_1c_2+c_3)^{\frac{N}{4}}
-c_1^{\frac{N}{2}-1}(c_1^2+c_2)^{\frac{N}{4}}\left[(r^2+1)c_1^2+\frac{N^2}{4}c_2
\right],$$
if $N\equiv 0 \mod 4$, while  if $N\equiv 2 \mod 4$ then we have
$$p_N(-rc_1)\equiv c_1^{\frac{N-2}{4}}(c^3_1 + c_2 c_1 + c_3)^{\frac{N+2}{4}}
- c_1^{\frac{N}{2}} (c_1^2+c_2)^{\frac{N-2}{4}}
\left((r^2+1) c_1^2+\frac{N^2+4}{4}c_2 \right).$$
In both cases, $p_N(-rc_1)\not \in I$.
\end{proposition}
\begin{proof}
The first part follows by substituting $t=-rc_1$ into the formulas in Corollary
\ref{pN-explicit} (and rearranging the terms), using the facts that $r$ is odd, that $2c_3=0$ and that
(see the proof of Proposition \ref{r-even})
$$2\cdot\left[-c_1\left(r^2c_1+\frac{N^2}{4}c_2 \right)\right]\in I.$$
To prove the last statement, it is enough to prove that $p_N(-rc_1)\not\in (I, c_2,c_3)$.
From the above formulas for $p_N(-rc_1)$, we get that (in both cases)
$$p_N(-rc_1)\equiv -r^2 c_1^{N+1}\mod (I,c_2,c_3).$$
On the other hand, from the definition of the ideal $I$, we get the inclusion
$$(I,c_2,c_3)=(2r(N-1)c_1, 2r^2c_1^2, c_2,c_3)\subset (2c_1,c_2,c_3).$$
Now we conclude by observing that if $r$ is odd then the element $-r^2 c_1^{N+1}$
does not belong to the ideal $(2c_1,c_2,c_3)$ and hence, a fortiori,
neither to the ideal $(I,c_2,c_3)$.
\end{proof}

\section{Some explicit computations}\label{explicit}

\subsection{Explicit generators}
\label{exp-gen}
In this section we give an explicit description for the generators of the Chow
ring $A^*(\HH)$, viewed as operational Chow ring (see \cite[Prop. 17, 19]{EGa}).
We need first to fix some notation and recall two auxiliary stacks introduced in \cite{AV}.
Consider a uniform ${\mathbb \mu_r}$-cover $\pi:\F\to S$ of a conic bundle
$p:\C\to S$ such that the ramification divisor
$W\subset\F$ and the branch divisor $D\subset \C$ are both finite and \'etale
over $S$ of degree $N$.
By the classical theory of cyclic covers and the Hurwitz formula, there exists
an $r$-root $\L^{-1}\in \Pic(\C)$ of $\O_{\C}(D)$
such that, called $f:\F\to \C$ the cyclic cover of degree $r$, we have
\begin{equation}\label{f2}
f^*(\L^{-1})=\O_{\F}(W),
\end{equation}
\begin{equation}\label{f3}
f_*(\O_{\F})=\bigoplus_{i=0}^{r-1}\L^i,
\end{equation}
\begin{equation}\label{f4}
\omega_{\F/S}=f^*(\omega_{\C/S})\otimes \O_{\F}((r-1)W).
\end{equation}
In the above notation, it is easy to check that $\HH$ is isomorphic to the stack
$\HHe$ whose fiber over a $k$-scheme $S$ is the groupoid of collections
$(\C\stackrel{p}{\to} S, \L, \L^{\otimes r}\stackrel{i}{\hookrightarrow}
\O_{\C})$, where the morphisms are natural Cartesian diagrams
(see \cite[section 2]{AV}). Moreover, we consider an auxiliary stack $\HHa$ whose fiber over a $k$-schemes $S$ is the groupoid of collections
$$\HHa(S)=\{(\C\stackrel{p}{\to} S, \L,
\L^{\otimes r}\stackrel{i}{\hookrightarrow} \O_{\C}, \phi:(\C,\L)\cong
(\P^1_S,\O_{\P^1_S} (-d)))\},
$$
where the isomorphism $\phi$ consists of an isomorphisms of $S$-schemes
$\phi_0:\C\cong \P^1_S$ plus an isomorphism of invertible sheaves
$\phi_1:\L\cong \phi_0^*\O_{\P^1_S}(-d)$.
In \cite[Theo. 4.1]{AV}, it is proved that $\HHa\cong \Asm$
and that the forgetful morphism $\HHa\to \HHe\cong \HH$ is a principal $GL_2/\mu_d$-bundle.
\begin{theorem}\label{gen-odd}
Assume that $d$ is odd. The Chow ring $A^*(\HH)$ is generated by the first two
Chern classes of the vector bundle of rank $2$
$$\E_1(\F\stackrel{\pi}{\to} S):=\pi_*
\omega_{\F/S}^{\frac{d-1}{2}}\left(\frac{1+r+d-N}{2}W \right).$$
\end{theorem}
\begin{proof}
The equivariant Chow ring $A^*_{GL_2/\mu_d}(\Asm)$ is a quotient of
the Chow ring of the group $GL_2/\mu_{d}$.
From the isomorphism of algebraic groups
$$
\begin{aligned}
GL_2/\mu_{d}&\stackrel{\cong}{\longrightarrow} GL_2\\
[A] &\longrightarrow (\det A)^{\frac{d-1}{2}} A,\\
\end{aligned}
$$
and the fact that $A^*_{GL_2}=\Z[c_1,c_2]$, where $c_1$ and $c_2$ are the Chern classes
of the standard representation of $GL_2$, we deduce that $A^*_{GL_2/\mu_{d}}(\Asm)$
is generated by the first two equivariant Chern classes of the vector bundle
$\tilde{\E_1}$ that associates to a trivial family $\P^1(V_S)=\P^1_S\stackrel{p_S}{\longrightarrow} S$ the vector bundle
$$\tilde{\E_1}(\P^1_S\stackrel{p_S}{\longrightarrow} S):=(\det V_S)^{\frac{d-1}{2}}\otimes V_S,
$$
where $V_S=V\times_k S$ and $V$ is the two dimensional standard $k$-representation of $GL_2$.
Clearly we have that $V_S=(p_S)_*(\O_{\P^1_S}(1))$. Moreover, from the Euler
exact sequence applied to the trivial family $p_S:\P^1_S\to S$
$$
0\to \O_{\P^1_S}\to p_S^*(V_S^*) (1)\to
\omega_{\P^1_S/ S}^{-1}\to 0,
$$
we deduce the $GL_2$-equivariant isomorphism
$p_S^*(\det V_S)\cong \O_{\P^1_S}(2)\otimes \omega_{\P^1_S/S}$.
Using the projection formula and the equality $(p_S)_*(\O_{\P^1_S})=\O_S$, we get
the $GL_2$-equivariant isomorphism
\begin{equation}\label{Euler}
\det V_S=(p_S)_*(\O_{\P^1_S}(2)\otimes \omega_{\P^1_S/S}).
\end{equation}
where we consider the canonical actions of $GL_2$ on $\P^1_S$ and
on the invertible sheaves involved. Using these two equalities and the isomorphism
$$
\phi:(\C,\L)\cong (\P^1_S,\O_{\P^1_S}(-d)),
$$
the $GL_2/\mu_{d}$-equivariant
vector bundle $\tilde{\E_1}$ on $\HHa$ descends on $\HHe$ to the vector bundle
$$
\E_1'(\C\stackrel{p}{\to} S, \L):= p_*\left(\omega_{\C/S}^{\frac{d-1}{2}}\otimes \L^{-1}\right).
$$
Consider now a $\mu_r$-cover $f:\F\to \C$ with ramification divisor $W$ as above. Using the formulas (\ref{f2}), (\ref{f4}), we have that $$f^*(\omega_{\C/S}^{(d-1)/2}\otimes \L^{-1})= \omega_{\F/S}^{(d-1)/2}\left(\frac{1+d+r-N}{2} W\right).$$ Therefore, using the projection formula and the formula (\ref{f3}), we get
$$
\pi_*\omega_{\F/S}^{\frac{d-1}{2}}\left(\frac{1+d+r-N}{2} W\right)=p_*\left( \left( \omega_{\C/S}^{\frac{d-1}{2}}\otimes \L^{-1}\right)\otimes \bigoplus_{i=0}^{r-1}\L^i  \right)=p_*\left(\omega_{\C/S}^{\frac{d-1}{2}}\otimes \L^{-1}\right),
$$
where in the last equality we used that $p_*(\omega_{\C/S}^{(d-1)/2}\otimes \L^j)=0$ for $j\geq 0$, which follows from the fact that $\omega_{\C/S}^{(d-1)/2}\otimes \L^j$ has negative degree on the fibers of $p$, if $j\geq 0$. This show that under the above isomorphism of stacks $\HHe\cong \HH$, the vector bundle $\E_1'$ goes into the vector bundle $\E_1$.
\end{proof}
\begin{theorem}\label{gen-even}
Assume that $d$ is even. The Chow ring $A^*(\H)$ is generated by the first Chern class of the line bundle
$$
\GG(\F\stackrel{\pi}{\to}S):=\pi_*\omega_{\F/S}^{d/2}\left(\frac{2+d-N}{2} W\right),
$$
and by the second and third Chern classes of the vector bundle of rank $3$
$$
\E_2(\F\stackrel{\pi}{\to} S):=\frac{\pi_* \omega_{\F/S}^{-1}((r-1)W)} {\pi_* \omega_{\F/S}^{-1}((r-2)W)}.
$$
Moreover $\pi_*\omega_{\F/S}^{-1}((r-2)W)=0$ if $d>2$.
\end{theorem}
\begin{proof}
The equivariant Chow ring $A^*_{GL_2/\mu_{d}}(\Asm)$ is a quotient of
the Chow ring of the group $GL_2/\mu_{d}$ which, for $d$ even, is isomorphic to
$$
\begin{aligned}
GL_2/\mu_{d}&\stackrel{\cong}{\longrightarrow} \G\times PGL_2\\
[A] &\longrightarrow ((\det A)^{d/2}, [A]).\\
\end{aligned}
$$

Recall that $A_{\G\times PGL_2}^*=\Z[c_1,c_2,c_3]/(2c_3)$, where $c_1$ is the first Chern class of the natural representation of $\G$ and $c_2, c_3$ are the second and third Chern classes of the representation ${\rm Sym}^2(V)$ of $PGL_2=SL_2/\mu_2$, where $V$ is the two dimensional standard $k$-representation of $SL_2$. Therefore we deduce that $A^*_{GL_2/\mu_{d}}(\Asm)$ is generated by the first Chern class of the line bundle $\tilde{\GG}$ that associates to a trivial family $\P^1(V_S)=\P^1(V\times_k S)=\P^1_S\stackrel{p_S}{\longrightarrow} S$ the $GL_2/\mu_{d}$-equivariant line bundle
$$\tilde{\GG}(\P^1_S\stackrel{p_S}{\longrightarrow} S):=(\det V_S)^{d/2},$$
and by the second and third Chern classes of the $GL_2/\mu_{d}$-equivariant vector bundle
$$\tilde{\E_2}(\P^1_S\stackrel{p_S}{\longrightarrow} S)={\rm Sym}^2(V_S).$$
Clearly we have that ${\rm Sym}^2(V_S)=(p_S)_*(\omega_{\P^1_S/S}^{-1})$. Moreover, from the formula (\ref{Euler}), we deduce that $(\det V_S)^{d/2}=(p_S)_*(\O_{\P^1_S}(d)\otimes \omega_{\P^1_S/S}^{d/2})$, where we consider the canonical actions of $GL_2/\mu_d$ on $\P^1_S$ and on the invertible sheaves involved. Using these two equalities and the isomorphism $\phi:(\C,\L)\cong (\P^1_S,\O_{\P^1_S}(-d))$, the $GL_2/\mu_{d}$-equivariant vector bundles $\tilde{\GG}$ and $\tilde{\E_2}$ on $\HHa$ descend on $\HHe$ to the vector bundles
$$
\GG'(\C\stackrel{p}{\to} S, \L):=
p_*(\omega_{\C/S}^{d/2}\otimes \L^{-1}),
$$
$$
\E_2'(\C\stackrel{p}{\to} S, \L):=p_*(\omega_{\C/S}^{-1}).
$$
Consider now the $\mu_r$-cover $f:\F\to \C$ with ramification divisor $W$ as above. Using the formulas (\ref{f2}), (\ref{f4}), we have that
$$
\begin{sis}
&f^*(\omega_{\C/S}^{d/2}\otimes \L^{-1})=
\omega_{\F/S}^{d/2}\left(\frac{2+d-N}{2}
W\right),\\
& f^*(\omega_{\C/S}^{-1})=\omega_{\F/S}^{-1}((r-1)W),\\
&f^*(\omega_{\C/S}^{-1}\otimes \L)=\omega_{\F/S}^{-1}((r-2)W).
\end{sis}
$$

Therefore, using the projection formula, the formula (\ref{f3}) and
the vanishings $p_*(\omega_{\C/S}^{d/2}\otimes \L^{-1}\otimes
\L^{i})=p_*(\omega_{\C/S}^{-1}\otimes
\L^{i+1})=0$ for $i\geq 1$ (because these line bundles have negative
degrees on the fibers of $p:\C\to S$), we get
$$\begin{sis}
&\pi_*\omega_{\F/S}^{d/2}\left(\frac{2+d-N}{2} W\right)=
p_*\left( \left( \omega_{\C/S}^{d/2}\otimes \L^{-1}\right)\otimes \bigoplus_{i=0}^{r-1}\L^i\right)=
p_*(\omega_{\C/S}^{d/2}\otimes \L^{-1}),\\
&\pi_*\omega_{\F/S}^{-1}\left((r-1) W\right)=
p_*\left(  \omega_{\C/S}^{-1}\otimes \bigoplus_{i=0}^{r-1} \L^i
\right)=p_*(\omega_{\C/S}^{-1})\oplus p_*(\omega_{\C/S}^{-1}\otimes \L),\\
&\pi_*\omega_{\F/S}^{-1}\left((r-2) W\right)=
p_*\left(\omega_{\C/S}^{-1}\otimes \L\otimes \bigoplus_{i=0}^{r-1} \L^i
\right)= p_*(\omega_{\C/S}^{-1}\otimes \L).\\
\end{sis}$$
This shows that, under the above isomorphism of stacks $\HHe\cong \HH$,
the vector bundles $\GG'$ and $\E_1'$ go respectively into $\GG$ and $\E_1$.
Moreover, if $d>2$ then $\omega_{\C/S}^{-1}\otimes \L$ has negative degree on
the fibers of $p:\C\to S$ and therefore we get the vanishing
$\pi_*\omega_{\F/S}^{-1}\left((r-2) W\right)=
p_*(\omega_{\C/S}^{-1}\otimes \L)=0.$
\end{proof}

\begin{remark}
\noindent
In the hyperelliptic case. i.e. for $r=2$, one recovers the results
of \cite[Theo. 4.1]{GV} for the Picard group and of \cite[Section 5.1]{EF}
for the Chow ring in the case $g$ even (i.e. $d=g+1$ odd).
\end{remark}

\subsection{$\lambda$-classes}

In this last part of the section, we want to express the tautological $\lambda$ classes
of Mumford (\cite{Mum})
in terms of the above explicit generators of $A^*(\HH)$. Recall that the lambda classes are
defined as Chern classes of the Hodge bundle:
$$\lambda_j(\F\stackrel{\pi}{\longrightarrow} S):=c_j(\pi_*(\omega_{\F/S})) \text{ for }
j=1, \cdots, g. $$

\begin{proposition}\label{lambda}
\noindent
\begin{itemize}
 \item[(i)] Assume that $d$ is odd.
Then the polynomial expressing $\lambda_j$ in terms of the above
generators $c_1(\E_1)$ and $c_2(\E_2)$ is the same as the one expressing
$$\bigoplus_{i=0}^{r-2} c_j\left((\det V)^{(r-1-i)\frac{1-d}{2}+1}
\otimes {\rm Sym}^{(r-1-i)d-2}(V)\right)$$
in terms of the generators $c_1$ and $c_2$ of $A^*_{GL_2}=\Z[c_1, c_2]$,
where $V$ is the standard two dimensional representation of $GL_2$.

\item[(ii)] Assume that $d$ is even. Then the polynomial expressing $\lambda_j$
in terms of the above generators $c_1(\GG)$, $c_2(\E_2)$ and $c_3(\E_2)$ is the same as the one expressing
$$\bigoplus_{i=0}^{r-2} c_j\left(W^{(r-1-i)d}\otimes
{\rm Sym}^{(r-1-i)d-2}(V)\right)$$
in terms of the generators $c_1$, $c_2$, and $c_3$ of $A^*_{\G\times PGL_2}=\Z[c_1, c_2, c_3]/(2 c_3)$,
where $W$ is the standard one dimensional representation of $\G$ and $V$ is the standard two dimensional
representation of $SL_2$.
\end{itemize}
\end{proposition}

\begin{proof}
Consider the Hodge vector bundle
$\Ho(\F\stackrel{\pi}{\longrightarrow} S):=\pi_*(\omega_{\F/S}).$
Using the formulas (\ref{f2}), (\ref{f3}), (\ref{f4}) and the standard
notations introduced in section \ref{exp-gen}, we compute
$$
\pi_*(\omega_{\F/S})=p_*f_* f^*(\omega_{\C/S}\otimes \L^{-(r-1)})=
p_*\left(\bigoplus_{i=0}^{r-1} \omega_{\C/S}\otimes \L^{i-(r-1)}\right)
=$$
$$=\bigoplus_{i=0}^{r-2} p_*\left(\omega_{\C/S}\otimes \L^{i-(r-1)}\right),
$$
where in the last equality we have used the vanishing $p_*(\omega_{\C/S})=0$.
For a trivial family $(\C \stackrel{p}{\longrightarrow} S, \L  )\cong
(\P^1_S\stackrel{p_S}{\longrightarrow} S,\: \O_{\P^1_S}(-d))$, using the formula
(\ref{Euler}), we get the $GL_2/\mu_{d}$-equivariant isomorphism
$$(p_S)_*(\omega_{\C/S}\otimes \L^{-1})=(p_S)_*\left( p_S^*(\det V_S)\otimes
\O_{\P^1_S}((r-1-i)d-2)\right)=$$
$$=\det V_S\otimes {\rm Sym}^{(r-1-i)d-2}(V_S),
$$
where $V_S=V\times_k S$ and $V$ is the standard representation of $GL_2$.

Suppose first that $d$ is odd. Using the isomorphism of algebraic groups
$$
\begin{aligned}
GL_2&\stackrel{\cong}{\longrightarrow} GL_2/\mu_{d}\\
A &\longrightarrow \left[(\det A)^{\frac{-d+1}{2d}} A\right],\\
\end{aligned}
$$
the above $GL_2/\mu_{d}$-representation $\det V_S\otimes
{\rm Sym}^{(r-1-i)d-2}(V_S)$ becomes isomorphic to the $GL_2$-representation
$$[(\det V)^{\otimes \frac{-d+1}{2d}2}\otimes (\det V) ]\otimes
[(\det V)^{\frac{-d+1}{2d}[(r-1-i)d-2]}\otimes {\rm Sym}^{(r-1-i)d-2}(V)]=$$
$$=(\det V)^{(r-i-1)\frac{1-d}{2}+1} {\rm Sym}^{(r-1-i)d-2}(V),$$ which gives the conclusion.

Finally, if $d$ is even then, using the isomorphism of algebraic groups
$$
\begin{aligned}
\G \times PGL_2&\stackrel{\cong}{\longrightarrow} GL_2/\mu_{d}\\
(\alpha, [A]) &\longrightarrow \alpha^{\frac{1}{d}} (\det A)^{-1/2} A,\\
\end{aligned}
$$
the above $GL_2/\mu_{d}$-representation $\det V_S\otimes
{\rm Sym}^{(r-1-i)d-2}(V_S)$
becomes isomorphic to the $\G\times PGL_2$-representation
$$[W^{\frac{2}{d}}]\otimes [W^{\frac{(r-1-i)d-2}{d}}  \otimes
{\rm Sym}^{(r-1-i)d-2}(V)]=W^{(r-1-i)d}\otimes
{\rm Sym}^{(r-1-i)d-2}(V), $$
which gives the conclusion.
\end{proof}

Consider now the natural representable map of stacks
$$\phi:\HH\to \Mg,$$
where $g=(r-1)(N-2)/2$.
It induces a pull-back map $\phi^*:\Pic(\Mg)\to \Pic(\HH)$. Recall
(see \cite{AC}) that if $g\geq 2$ then $\Pic(\Mg)$ is cyclic generated
(freely if $g\geq 3$) by $\lambda_1$.

\begin{corollary}\label{lambda_1}
The class $\lambda_1$ in $\Pic(\HH)=\langle c_1\rangle$ is equal to
$$\lambda_1=\begin{cases}
\sum_{j=0}^{r-1} \binom{dj}{2} c_1 &\text{ if } d \text{ is odd,}\\
2\sum_{j=0}^{r-1} \binom{dj}{2} c_1 &\text{ if } d \text{ is even.}\\
\end{cases}$$
\end{corollary}
\begin{proof}
Assume first that $d$ is odd. Consider the formula for $\lambda_1$ given
in Proposition \ref{lambda}(i). Using the relations
$c_1(\det(V)^m)=m c_1(\det(V))=m c_1$ and $c_1(\Sym^m(V))=\frac{(m+1)m}{2}
c_1(V)=\frac{(m+1)m}{2} c_1$, we get
\begin{eqnarray*}
\lambda_1&=&\sum_{i=0}^{r-2} \left\{\left[(r-1-i)d-1\right]\left[(r-1-i)\frac{1-d}{2}+1\right] +
\right. \\
&+& \left. \frac{[(r-1-i)d-2][(r-1-i)d-1]}{2} \right\}c_1\\
&=& \sum_{i=0}^{r-2}\binom{(r-1-i)d}{2}c_1=\sum_{j=0}^{r-1}
\binom{jd}{2}c_1.
\end{eqnarray*}

Assume now that $d$ is even. Consider the formula for $\lambda_1$ given
in Proposition \ref{lambda}(ii). Using the relations
$c_1(W^m)=m c_1(W)=m c_1$ and $c_1(\Sym^m(V))=0$, we get
$$\lambda_1=\sum_{i=0}^{r-2} [(r-1-i)d-1][(r-1-i)d]c_1=
2\sum_{j=0}^{r-1} \binom{jd}{2}c_1.$$
\end{proof}

\begin{remark}
In the hyperelliptic case $\H=\mathcal{H}_{\rm sm}(1,2,g+1)$,
one recovers the result of \cite[Corollary 4.4]{GV} since
$$\lambda_1=
\begin{cases}
\binom{g+1}{2}c_1\equiv \frac{g}{2} c_1 \mod 2(2g+1)
& \text{ if } g \text{ is even, }\\
2\binom{g+1}{2}c_1\equiv g \,c_1 \mod 4(2g+1)& \text{ if } g \text{ is odd. }
\end{cases}
$$
\end{remark}

\section*{Acknowledgements} We thank Dan Edidin, Sergey Gorchinskiy, Burt Totaro and Angelo Vistoli for helpful discussions. We are grateful to the referee
for useful suggestions.

\end{document}